\documentclass{elsarticle}

\usepackage{amssymb}
\usepackage{amsthm}

\usepackage{amsmath}
\usepackage{subcaption}
\usepackage{listings}
\usepackage{hyperref}
\usepackage{multirow}
\graphicspath{ {./Figures/} }

\begin{document}
	
	\begin{frontmatter}
		
		
		
		\title{Explicit Runge-Kutta-Chebyshev methods of second order \\ with monotonic stability polynomial}
		\author[bsu]{Boris Faleichik}
		\ead{faleichik@bsu.by}
		\cortext[corauthor]{Corresponding author}
		
		\author[sunbim]{Andrew Moisa}
		\ead{andrey.moysa@gmail.com}
		
		\affiliation[bsu]{organization={Belarusian State University},
			addressline={Nezavisimosti ave., 4}, 
			city={Minsk},
			postcode={220030}, 
			country={Belarus}}
		\affiliation[sunbim]{organization={Sunbim Ltd.},
			addressline={71-75 Shelton Street, Covent Garden}, 
			city={London},
			postcode={WC2H~9JQ}, 
			country={United~Kingdom}}

		\begin{abstract}
			A new Chebyshev-type family of stabilized explicit methods for solving mildly stiff ODEs is presented. Besides conventional conditions of order and stability we impose an additional restriction on the methods: their stability function must be monotonically increasing and positive along the largest possible interval of negative real axis. Although stability intervals of the proposed methods are smaller than those of classic Chebyshev-type methods, their stability functions are more consistent with the exponent, they have more convex stability regions and smaller error constants. These properties allow the monotonic methods to be competitive with contemporary stabilized second-order methods, as the presented results of numerical experiments demonstrate.   
		\end{abstract}
		
		\begin{graphicalabstract}
			\begin{center}
				\includegraphics{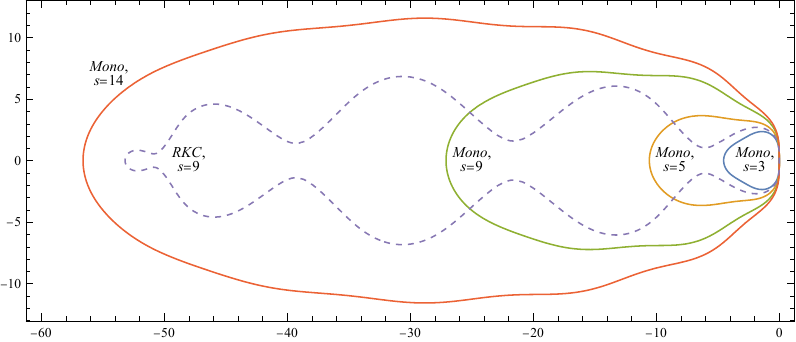}
			\end{center}		
		\end{graphicalabstract}
		
		
		
	\end{frontmatter}
	
	
	
	\section{Introduction}
	We assume that the reader is familiar with basic concepts of stability analysis for numerical ODE solution methods \cite[IV.2]{hairer2} and techniques of constructing stabilized explicit Runge--Kutta methods \cite[IV.2]{hairer2}, \cite[V]{hundsdorfer}, \cite{rkc_conv}, \cite{rkc}.
	
	As is customary when building stabilized RK methods for systems of ordinary differential equations
	\begin{equation}
		y'(x)=f(x, y(x)),\quad y(x_0) = y_0,
	\end{equation}
	we start from constructing stability polynomial $R_s$ of degree $s$ with required properties. The key condition is 
	\begin{equation}\label{cond:mono}
		R'_s(x) \geq 0\quad \text{and}\quad R_s(x) > 0\quad \forall x\in (-\rho_s,0], 
	\end{equation}
	where $\rho_s$, which will be called the length of monotonicity interval, should be maximized. The motivation for this condition goes back to \cite{bobkov}: we pursuit improved consistency between the stability function and the exponent, since in the case of a linear system ${y'(x) = A y(x)}$ exact and approximate solutions are related as
	\begin{equation*}
		y(x_0 + h) = \exp(h A) y_0 \approx R_s(hA)y_0 = y_1.
	\end{equation*}
	We also impose order conditions 
	\begin{equation}\label{cond:order}
		R_s(0) = R'_s(0) = R''_s(0) = 1,
	\end{equation}
	which is sufficient for the final Runge-Kutta scheme to have order two. The following layout is similar to the one used in deriving RKC methods \cite{rkc_conv}.

	\section{Method construction}

	Maximizing $\rho_s$ implies that, besides non-negativity, $R'_s$ should also have the minimum possible deviation from zero within the interval of monotonicity $[-\rho_s, 0]$. Thus, it is natural to shape $R'_s$ in the form of a shifted and scaled Chebyshev polynomial of the first kind:
	\begin{equation}\label{dRs}
		R'_s(x)= b_{s-1} \bigl(1 + T_{s-1}(w_0+w_1 x)\bigr), 
	\end{equation}
	where 
	\begin{equation}
		b_s = \frac{1}{1+T_s(w_0)},
	\end{equation}
	which directly follows from the order condition $R'_s(0) = 1$. Now we can define $\rho_s$ from the scaling identity $w_0-w_1 \rho_s = -1$:
	\begin{equation}
	 	\rho_s = \frac{1+w_0}{w_1},
	\end{equation} 
	and apply the remaining order conditions to find the parameters $w_0$ and $w_1$. Condition $R''(0)=1$ gives
	\begin{equation}\label{o2cond}
		w_1=\frac{1+T_{s-1}(w_0)}{T'_{s-1}(w_0)} = \frac{1}{b_{s-1} T_{s-1}'(w_0)}.
	\end{equation}
	From \eqref{dRs} we put 
	\begin{equation}
		R_s(x) = b_{s-1} \int_{-\rho_s}^x\bigl(1+T_{s-1}(\xi)\bigr)d\xi  
	\end{equation}
	and utilize one of the well-known properties of Chebyshev polynomials
	\begin{equation*}
		\int T_s(x)dx = \frac{1}{2(s+1)} T_{s+1}(x) - \frac{1}{2(s-1)} T_{s-1}(x).
	\end{equation*}
	The result is
	\begin{equation*}
		R_s(x)=\alpha_s+b_{s-1} x + \gamma_s T_s(w_0+w_1 x) + \delta_s T_{s-2}(w_0+w_1 x),
	\end{equation*}
	where
	\begin{equation}\label{agd}
		\alpha_s = \frac{b_{s-1}}{w_1}\left(1 + \frac{(-1)^s }{s(s-2) }+ w_0
		\right) ,\quad
		\gamma_s = \frac{b_{s-1}}{2s w_1},\quad \delta_s = -\frac{b_{s-1}}{2(s-2)w_1}.
	\end{equation}
	From $R(0)=1$ it follows that 
	\begin{equation}\label{o0cond}
		\alpha_s = 1 - \gamma_s T_s(w_0) - \delta_s T_{s-2}(w_0),
	\end{equation} 
	and 
	\begin{equation}\label{Rs1}
		R_s(x) = 1 + b_{s-1} x + \gamma_s \bigl(T_s(w_0 + w_1 x) - T_s(w_0)\bigr) + 
		\delta_s  \bigl(T_{s-2}(w_0 + w_1 x) - T_{s-2}(w_0)\bigr),
	\end{equation}
	Eliminating $w_1$ from \eqref{agd}, \eqref{o0cond} and using \eqref{o2cond} we get the last equation which allows to determine $w_0$: 
	\begin{equation}\label{w0}
		1+ \frac{(-1)^s}{s(s-2)}+w_0+\frac{1}{2s}T_s(w_0)-\frac{1}{2(s-2)}T_{s-2}(w_0) = \frac{(1+T_{s-1}(w_0))^2}{T_{s-1}'(w_0)}.
	\end{equation}
	In practice we solve this equation numerically in Wolfram Language using high-precision arithmetic with 500 significant decimal places.
	
	Just like most of the existing Runge-Kutta-Chebyshev methods, the numerical RK scheme implementing the stability polynomial we have built will be based on the basic three-term recurrence relation 
	\begin{equation*}
		T_s(x) = 2x T_{s-1}(x) - T_{s-2}(x). 	
	\end{equation*} 
	To use this formula in \eqref{Rs1} and maintain internal stability of the corresponding RK stages let us introduce the polynomials
	\begin{equation}
		\tilde R_j(x) = 1 + b_j\bigl(T_j(w_0+w_1 x) - T_j(w_0)\bigr),\quad j = 0, 1, \ldots, s.
	\end{equation}
	Then $R_s$ takes the form 
	\begin{equation}\label{Rs2}
		R_s(x) = 1 + b_{s-1} x + \frac{\gamma_s}{b_s} \bigl(\tilde R_s(x) - 1\bigr) + 
		\frac{\delta_s}{b_{s-2}}  \bigl(\tilde R_{s-2}(x) - 1\bigr),
	\end{equation}
	and the three-term recurrence gives
	\begin{equation}
	 	\tilde R_j(x) = 1 + \mu_j \bigl(\tilde R_{j-1}(x) - 1\bigr) + \nu_j \bigl(\tilde R_{j-2}(x) - 1\bigr) + 
	 	\tilde \mu_j x \bigl(\tilde R_{j-1}(x) - b_{j-1}\bigr),
	 	\quad j \geq 2,
	\end{equation} 
	where
	\begin{equation}\label{munu}
		\mu_j = 2w_0\frac{b_j}{b_{j-1}},\quad \nu_j = -\frac{b_j}{b_{j-2}},\quad \tilde\mu_j = 2w_1\frac{b_j}{b_{j-1}}.
	\end{equation}

\begin{figure}
		\caption{Stability polynomials for ${s=9}$. Left: comparison of monotonic and RKC method. Right: graphs of the internal stability polynomials $\tilde R_j$~(colored) and $R_9$~(gray).}
		\label{fig1}
		\begin{center}
			\includegraphics{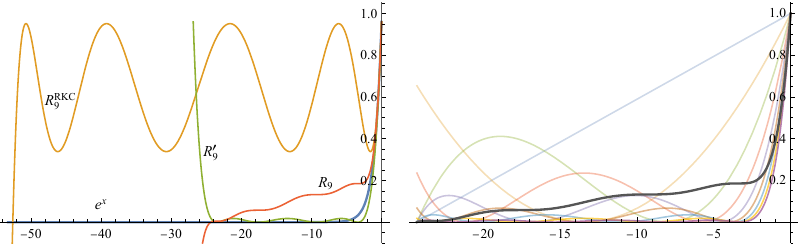}
		\end{center}
	\end{figure}

	The resulting numerical scheme is very similar to the second-order RKC method \cite{rkc_conv}:
	\begin{subequations}\label{method}
			\begin{align}
			&Y_0 = y_0,\\
			&Y_1 = y_0 + h\, b_1 w_1 F_0,\\
			&Y_j = (1-\mu_j-\nu_j) y_0 + \mu_j Y_{j-1} + \nu_j Y_{j-2} + h\,\tilde\mu_j(F_{j-1} - b_{j-1} F_0),\quad j = 2,\ldots,s,\\
			&y_1 = \left(1-\frac{\gamma_s}{b_s} - \frac{\delta_s}{b_{s-2}}\right)y_0 + \frac{\gamma_s}{b_s} Y_s + \frac{\delta_s}{b_{s-2}}Y_{s-2} + h\, b_{s-1} F_0,
		\end{align}
	\end{subequations}	where $F_i = f(x_0+c_i h, Y_i)$, and
	\begin{equation}\label{cj_1}
		c_0=0,\quad c_1 = w_1 b_1,\quad c_j = \mu_j c_{j-1} + \nu_j c_{j-2} + \tilde\mu_j(1-b_{j-1}), \quad j = 2, \ldots, s-1.
	\end{equation}
	Stability functions for each stage value $Y_j$ are equal to $\tilde R_j$. It is easy to check that 
	\begin{equation}
		R_s'(x)=\tilde R_{s-1}(x),
	\end{equation}
	so by design $Y_{s-1}$ is a stable first-order approximation of $y(x_0+h)$ and can be used for error estimation (although in practice we use another approach, see below). 
	The internal stability of the method is conditioned by the fact that, by construction, $\tilde R_j$ satisfy condition 
	$0 \leq \tilde R_j(x) < 1$ $\forall x\in[-\rho_s,0)$, $j > 0$. From \eqref{cj_1} we also have 
	\begin{equation}
		c_j = w_1 b_j T_j'(w_0),
	\end{equation}
	which leads to
	$$c_{j-1} < c_{j} < c_{s-1} = 1,\quad \forall j = 1,\ldots, s-2.$$
	Graphs of the final and the internal stability functions for $s=9$ are shown in Figure \ref{fig1}.

	\section{Properties of monotonic methods}

\begin{figure}
		\caption{Comparison of stability regions for second-order monotonic and RKC methods.}\label{fig3}
		\begin{center}
			\includegraphics{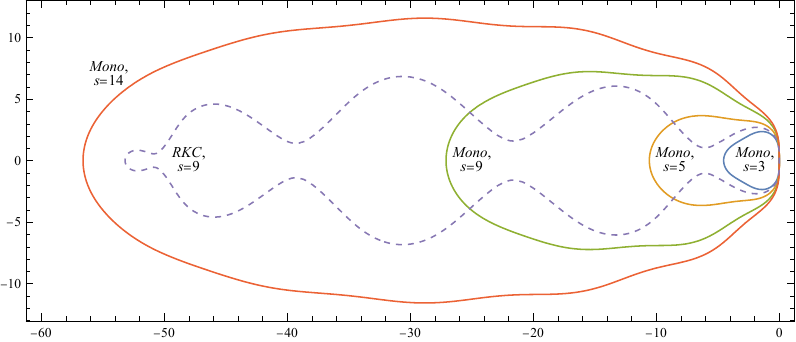}
		\end{center}

	\end{figure}

	\paragraph{Monotonicity interval} Due to the complexity of equation \eqref{w0} we do not have an exact formula for $\rho_s$. A numerical fit for a set of calculated methods gives  
	\begin{equation*}
		\rho_s \approx 0.31 \cdot (s+0.83)^{1.87},
	\end{equation*}
	see \eqref{s_est} below. Exact values for selected $s$ are shown in Table~\ref{table1}. As for stability intervals, it is clear that by construction their length for large $s$ is just slightly greater than $\rho_s$. Examples of stability regions and comparison with second-order RKC method are shown in Figure~\ref{fig3}. 

\begin{table}
	\caption{Parameters of some monotonic methods}
	{\small
	\begin{equation*}
		\begin{array}{r|lllllll} 
		s & \rho_s & C_s & w_0 & w_1 & b_{s-1} & \gamma_s & -\delta_s  \\ \hline
  3    & 3.5874010 & 0.0833333 & 1.2599210 & 0.62996052 & 0.31498026 & 0.08333333 & 0.25         		   \\
  5    & 8.6189019 & 0.0510313 & 1.4915378 & 0.28907833 & 0.04202332 & 0.01453700 & 0.02422833 		   \\
  10   & 29.268039 & 0.0322256 & 1.2057371 & 0.07536333 & 0.00679083 & 0.00450539 & 0.00563174 		   \\
  20   & 100.80657 & 0.0239240 & 1.0734470 & 0.02056856 & 0.00143509 & 0.00174428 & 0.00193809 		   \\
  50   & 525.59171 & 0.0183733 & 1.0175279 & 0.00383858 & 0.00021006 & 0.00054724 & 0.00057004 		   \\
  100  & 1855.5228 & 0.0158146 & 1.0057090 & 0.00108094 & 0.00005116 & 0.00023664 & 0.00024147 		   \\
  200  & 6617.5217 & 0.0139362 & 1.0018102 & 0.00030250 & 0.00001263 & 0.00010444 & 0.00010549 		   \\
  500  & 36059.771 & 0.0120702 & 1.0003830 & 0.00005547 & 2.008\cdot 10^{-6} & 0.00003620 & 0.00003634 \\
  1000 & 131320.58 & 0.0109659 & 1.0001157 & 0.00001523 & 5.010\cdot 10^{-7} & 0.00001644 & 0.00001648 \\
  2000 & 481823.56 & 0.0100482 & 1.0000344 & 4.150\cdot 10^{-6} & 1.251\cdot 10^{-7} & 7.536\cdot 10^{-6} & 7.543\cdot 10^{-6}  \\
		\end{array}
	\end{equation*}
	}
	\label{table1}
	\end{table}

	\paragraph{Error constant of stability function}

\begin{figure}
		\caption{Monotonic vs RKC methods. Left: lengths of stability intervals. Right: values of 
		$R^{(3)}_s(0)$.}\label{fig2}
		\begin{center}
			\includegraphics{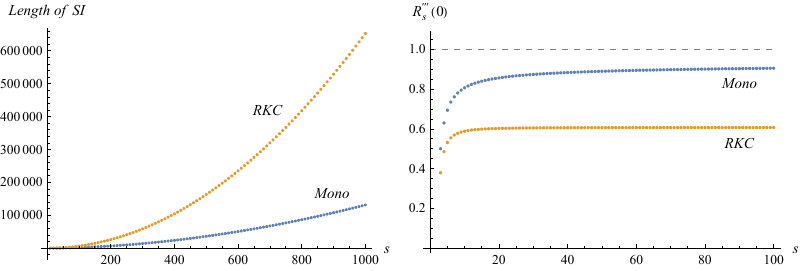}
		\end{center}
	\end{figure}

	We define the error constant of stability function as $C_s = \bigl(1-R_s'''(0)\bigr)/6$. From Figure~\ref{fig2} we can see that the error constants of monotonic methods are much smaller than those of RKC methods. For large $s$ the difference is about 6 times, see also Table~\ref{table1}.

	\section{Implementation details}


	To get a formula for estimating the required number of stages we took a set of numerically calculated points $(\rho_s, s)$ and fitted a nonlinear model of the form $s = a + b \rho ^c$. The result is 
	\begin{equation}\label{s_est}
		s \approx -0.8306782178712795 +  1.8547887825836553\cdot \rho^{0.533871357807877},	
	\end{equation}
	where $\rho$ is the spectral radius estimation.
	
	The first way to estimate the local error was already mentioned: the value of $Y_{s-1}$ from \eqref{method} can be used as an embedded method. However, its error constant varies when $s$ changes. Therefore, to reduce the number of rejected steps when changing $s$ some "normalization" is needed \cite{tsrkc2}. We will not delve into this process because another approach has been chosen in our current implementation.
	
	The other way to estimate the local error is based on a Taylor series expansion of the local solution \cite{rkc}, \cite{tsrkc3}. Let us suppose there is a first order embedded method which has the leading term of the local error expansion of the form
	\begin{equation}
		le = \frac{1}{10} h^2 \frac{d^2 y(x_0)}{dx^2}.
	\end{equation} 
	Then, an asymptotically correct estimate
	\begin{equation}\label{errest}
		err = \frac{1}{10} \left( y_0 - y_1 + f(x_0 + h, y_1) \right)
	\end{equation}
	is used for the practical error estimation. All other details of the step size selection are borrowed from \cite{rkc}.

	\section{Numerical experiments}
	\label{numexps}
	
	We conclude this paper by presenting some results of numerical experiments. The monotonic methods are implemented in a solver called MONO in C programming language. It is compared to RKC solver \cite{rkc} and the current version of TSRKC2 methods \cite{tsrkc2}, \cite{moisa_faleichik}, \cite{stabmethods}. Source code and all the examples are available at the repository \url{https://github.com/MoisaAndrew/StabilizedMethods}. Experiments with the code are welcome.
	
	We chose the following stiff problems:
	\begin{enumerate}
		\item CUSP: a combination of Zeeman's "cusp catastrophe" model combined with the 
		van der Pol oscillator \cite[pp. 147-148]{hairer2}.
		\item FINAG: The FitzHug and Nagumo nerve conduction equation, converted
		into ODE’s by the method of lines \cite{rock4}.
		\item BURGERS: Bateman–Burgers equation \cite{rock4}.
		\item COMB: a scalar two-dimensional nonlinear (hotspot) problem from combustion theory \cite[p. 439]{hundsdorfer}.
	\end{enumerate}
	
	All parameters and output points for these problems are directly borrowed from the source papers. The results are presented in Table~\ref{tab:stats}. It can be seen that the numerical properties of the proposed method coincide with its theoretical properties: it uses larger step sizes to obtain the same accuracy, however, a greater number of internal stages is required to obtain the same stability interval length. In general, it behaves similarly to other second order stabilized methods, though, in case of COMB problem, monotonic method has a better match between required and actual accuracy.
	
	\begin{table}
		\centering
		\caption{Results  for the compared methods; $err$ is Euclidean norm global error, $N_f$ is total number of function evaluations, $N_{accpt}$ is number of accepted steps, $N_{rejct}$ is number of rejected steps and $T_{comp}$ is computing time in milliseconds}
		\begin{tabular}[c]{|l|l|l|lrrrr|}
			\hline
			Problem & Method & tol & $err$ & $N_{f}$ & $N_{accpt}$ & $N_{rejct}$ & $T_{comp}$ \\
			\hline
			\hline
			\multirow{9}{*}{CUSP} & \multirow{3}{*}{RKC} & $10^{-3}$ & $4.61 \cdot 10^{-3}$ & 1811 & 83 & 2 & 0.18 \\
			& & $10^{-5}$ & $2.37 \cdot 10^{-4}$ & 3765 & 309 & 4 & 0.40 \\
			& & $10^{-7}$ & $1.14 \cdot 10^{-5}$ & 8640 & 1303 & 5 & 1.00 \\
			\cline{2-8}
			& \multirow{3}{*}{TSRKC2} & $10^{-3}$ & $1.97 \cdot 10^{-4}$ & 4584 & 262 & 53 & 0.44 \\
			& & $10^{-5}$ & $4.30 \cdot 10^{-6}$ & 9163 & 1289 & 95 & 1.02 \\
			& & $10^{-7}$ & $3.87 \cdot 10^{-7}$ & 26117 & 8187 & 221 & 3.44 \\
			\cline{2-8}
			& \multirow{3}{*}{MONO} & $10^{-3}$ & $2.42 \cdot 10^{-4}$ & 3809 & 71 & 3 & 0.37 \\
			& & $10^{-5}$ & $1.31 \cdot 10^{-5}$ & 8494 & 427 & 2 & 0.85 \\
			& & $10^{-7}$ & $5.14 \cdot 10^{-7}$ & 24420 & 3796 & 5 & 2.87 \\
			\hline
			\hline
			\multirow{8}{*}{FINAG} & \multirow{2}{*}{RKC} & $10^{-5}$ & $5.36 \cdot 10^{-1}$ & 2617 & 340 & 7 & 0.88 \\
			& & $10^{-7}$ & $2.62 \cdot 10^{-2}$ & 6631 & 1716 & 0 & 2.61 \\
			\cline{2-8}
			& \multirow{3}{*}{TSRKC2} & $10^{-3}$ & $5.75 \cdot 10^{+0}$ & 1801 & 205 & 7 & 0.54 \\
			& & $10^{-5}$ & $8.61 \cdot 10^{-2}$ & 4890 & 1283 & 0 & 1.74 \\
			& & $10^{-7}$ & $1.68 \cdot 10^{-3}$ & 19738 & 8869 & 0 & 8.83 \\
			\cline{2-8}
			& \multirow{3}{*}{MONO} & $10^{-3}$ & $4.50 \cdot 10^{+0}$ & 2673 & 76 & 14 & 0.81 \\
			& & $10^{-5}$ & $1.21 \cdot 10^{-1}$ & 4654 & 505 & 2 & 1.59 \\
			& & $10^{-7}$ & $2.61 \cdot 10^{-3}$ & 17413 & 4800 & 0 & 7.34 \\
			\hline
			\hline
			\multirow{9}{*}{BURGERS} & \multirow{3}{*}{RKC} & $10^{-3}$ & $3.41 \cdot 10^{-2}$ & 277 & 46 & 7 & 0.12 \\
			& & $10^{-5}$ & $1.95 \cdot 10^{-3}$ & 466 & 118 & 4 & 0.21\\
			& & $10^{-7}$ & $1.52 \cdot 10^{-4}$ & 1094 & 462 & 0 & 0.58 \\
			\cline{2-8}
			& \multirow{3}{*}{TSRKC2} & $10^{-3}$ & $4.80 \cdot 10^{-2}$ & 289 & 99 & 1 & 0.13 \\
			& & $10^{-5}$ & $6.93 \cdot 10^{-4}$ & 573 & 242 & 0 & 0.29 \\
			& & $10^{-7}$ & $9.82 \cdot 10^{-6}$ & 3920 & 1959 & 0 & 2.14 \\
			\cline{2-8}
			& \multirow{3}{*}{MONO} & $10^{-3}$ & $3.84 \cdot 10^{-2}$ & 265 & 20 & 2 & 0.10 \\
			& & $10^{-5}$ & $1.17 \cdot 10^{-3}$ & 505 & 109 & 0 & 0.23 \\
			& & $10^{-7}$ & $1.75 \cdot 10^{-5}$ & 3224 & 1074 & 0 & 1.54 \\
			\hline
			\hline
			\multirow{9}{*}{COMB} & \multirow{3}{*}{RKC} & $10^{-3}$ & $1.84 \cdot 10^{+1}$ & 979 & 62 & 2 & 12.39 \\
			& & $10^{-5}$ & $1.20 \cdot 10^{+0}$ & 1954 & 290 & 0 & 25.59 \\
			& & $10^{-7}$ & $5.97 \cdot 10^{-2}$ & 4745 & 1491 & 0 & 65.04 \\
			\cline{2-8}
			& \multirow{3}{*}{TSRKC2} & $10^{-3}$ & $1.84 \cdot 10^{+1}$ & 764 & 63 & 1 & 9.14 \\
			& & $10^{-5}$ & $3.77 \cdot 10^{-1}$ & 2599 & 651 & 0 & 34.35 \\
			& & $10^{-7}$ & $4.17 \cdot 10^{-3}$ & 14997 & 6510 & 0 & 224.53 \\
			\cline{2-8}
			& \multirow{3}{*}{MONO} & $10^{-3}$ & $3.72 \cdot 10^{-1}$ & 2167 & 39 & 7 & 26.31 \\
			& & $10^{-5}$ & $1.81 \cdot 10^{-2}$ & 2975 & 355 & 0 & 39.28 \\
			& & $10^{-7}$ & $6.12 \cdot 10^{-4}$ & 13993 & 3563 & 0 & 196.54 \\
			\hline
		\end{tabular}
		\label{tab:stats}
	\end{table}


	\bibliographystyle{elsarticle-num-names} 
	\bibliography{references}
	
	
	
	
	
\end{document}